\documentclass{article}
\usepackage{latexsym}
\usepackage{amssymb}
\usepackage{graphicx}
\newtheorem{thm}{Theorem}

\newtheorem{prp}[thm]{Proposition}

\begin{document}

\title{Cospectral regular graphs with \\ and without a perfect matching}
\author{Zolt\'an L. Bl\'azsik, \\
{\small Department of Computer Science,}\\
{\small E\"{o}tv\"{o}s University, Budapest,}\\
%Pázmány P. s. 1/c, Budapest, H-1117, Hungary
{\small\tt blazsik@cs.elte.hu}
\\[5pt]
Jay Cummings, \\
{\small Department of Mathematics,} \\
{\small University of California, San Diego, USA,}\\
% (UCSD)
%9500 Gilman Drive # 0112
%La Jolla, CA  92093-0112
%United States of America (USA)
{\small\tt jjcummings@math.ucsd.edu}
\\[3pt]
and
\\[3pt]
Willem H. Haemers
\\
{\small Department of Econometrics and O.R.,}
\\
{\small Tilburg University, Tilburg, The Netherlands,}
\\
{\small\tt haemers@uvt.nl}
}
\date{}
\maketitle
\abstract{
\noindent
For each $b\geq 5$ we construct a pair of cospectral $b$-regular graphs,
where one has a perfect matching and the other one not.
This solves a research problem posed by the third author at the 22nd British Combinatorial Conference.
}
\\[3pt]
\\[5pt]

\section{Introduction}
By K\H{o}nig's theorem, regular bipartite graphs of positive degree have a perfect matching.
For regular graphs which are not bipartite there exists a powerful sufficient condition for existence
of a perfect matching in terms of the spectrum of the adjacency matrix; see \cite{BH,CGH,H}.
Bipartiteness as well as regularity can be deduced from the spectrum (see~\cite{vDH}).
So it seems natural to ask whether for a regular graph existence of a perfect matching can be seen from the spectrum.
At the 22nd British Combinatorial Conference the third author (believing that the answer should be negative),
posed the problem of finding two cospectral regular graph, one with a perfect matching and one without one 
(see~\cite{C}, Problem 22.8).

For non-regular graphs there exist easy examples.
The disjoint union of the $4$-cycle $C_4$ and the path $P_{n-4}$ has a perfect matching when $n$ is even,
and is cospectral with a graph consisting of the path $P_{n-4}$ with two pendant vertices attached to each endpoint,
which obviously has no perfect matching.
More interesting connected examples where found by Aiden Roy (private communication).

Below we will construct a pair of connected $b$-regular graphs where one has a perfect matching,
and the other one not, for every $b\geq 5$.
The smallest example is a pair of $5$-regular graphs on $42$ vertices.
%The $6$-regular example contains $44$ vertices.
In general, if $b$ is odd the example below contains $b^2+5b-8$ vertices,
and if $b$ is even the example contains $b^2+3b-10$ vertices.

The main tool is the following result of Godsil and McKay~\cite{GM} (see also~\cite{vDH}).
\begin{prp}
Let $G$ be a graph and let $\{X,Y\}$ be a partition of the vertex set.
Suppose that $X$ induces a regular subgraph, and that each vertex $y \in Y$ has $0, |X|/2,$ or $|X|$ neighbors in $X$.
Make a new graph $G'$ from $G$ as follows.
For each $y\in Y$ with $|X|/2$ neighbors in $X$, delete the $|X|/2$ edges between $y$ and $X$, 
and join $y$ to the $|X|/2$ other vertices of $X$.
Then $G$ and $G'$ are cospectral.
\end{prp}

The set $X$ is called a \emph{switching set}.
The operation that changes $G$ to $G'$ is called \emph{Godsil-McKay switching}.

\section{Construction}

\begin{thm}
For each $b\geq 5$ there exists a pair of cospectral connected $b$-regular graphs,
where one has a perfect matching and the other one not.
\end{thm}
{\bf  Proof.}
We will prove the theorem by constructing a $b$-regular graph with a Godsil-McKay switching set $X$ and
no perfect matching, for which switching will introduce many perfect matchings.

First assume $b$ is odd.
Define the graph $H_b$ to be the complement of the disjoint union of $(b-1)/2$ edges and the path $P_2$.
Then $H_b$ has $b+2$ vertices, and each vertex has degree $b$ except for one vertex $u$ of degree $b-1$.
To $u$ we attach a pendant edge $\{u,v\}$, which increases its degree to $b$.
Call the graph thus obtained $\widetilde{H_b}$.
Notice that $H_b$ has an odd number of vertices, and therefore no perfect matching, while $\widetilde{H_b}$
has many perfect matchings, each of which contains the edge $\{u,v\}$.
Consequently, by attaching any other graph $F$ to $\widetilde{H_b}$ by identifying $v$ with some vertex in $F$,
the result has the property that no edge in $F$ which is incident with $v$ can be in a perfect matching.

We define the graph on the switching set $X$ to be $K_3+C_{2b-5}$,
the disjoint union of a triangle and a cycle with $2b-5$ vertices.
The construction of $Y$ starts with $b-2$ disjoint copies of $\widetilde{H_b}$.
We define $W$ to be the set of vertices consisting of the $b-2$ copies of $v$.
Each $w\in W$  will be joint to $b-1$ vertices of $X$, such that $w$ is joint to every vertex of the
triangle and no two vertex degrees of the larger cycle differ by more than one.
Notice that our graph is now connected, every vertex except those in the larger cycle has degree $b$,
and every vertex in $Y$ is adjacent to $0$, $|X|/2$, or $|X|$ vertices in $X$, so $X$ is a switching set.

We will enlarge $Y$ and add $(b-2)(b-1)$ edges between $Y$ and $X$ such that $X$ remains a switching set, and
each vertex gets degree $b$, as desired.
To this end, first add one more copy of $\widetilde{H_b}$ and insert $b-1$ edges between the copy of $v$ and the
vertices in $X$ belonging to the larger cycle, such that the degrees of these vertices still differ by at most one.
Next we add $(b-3)/2$ disjoint edges, and join both endpoints of each of these edges to $b-1$ vertices of the larger cycle in $X$, 
such that the degree of the vertices in $X$ become equal to $b$.
The result is shown Figure~1.
\begin{figure}[h]
\vspace{-20pt}
\hspace{-20pt}
\includegraphics[width=450pt,height=600pt]{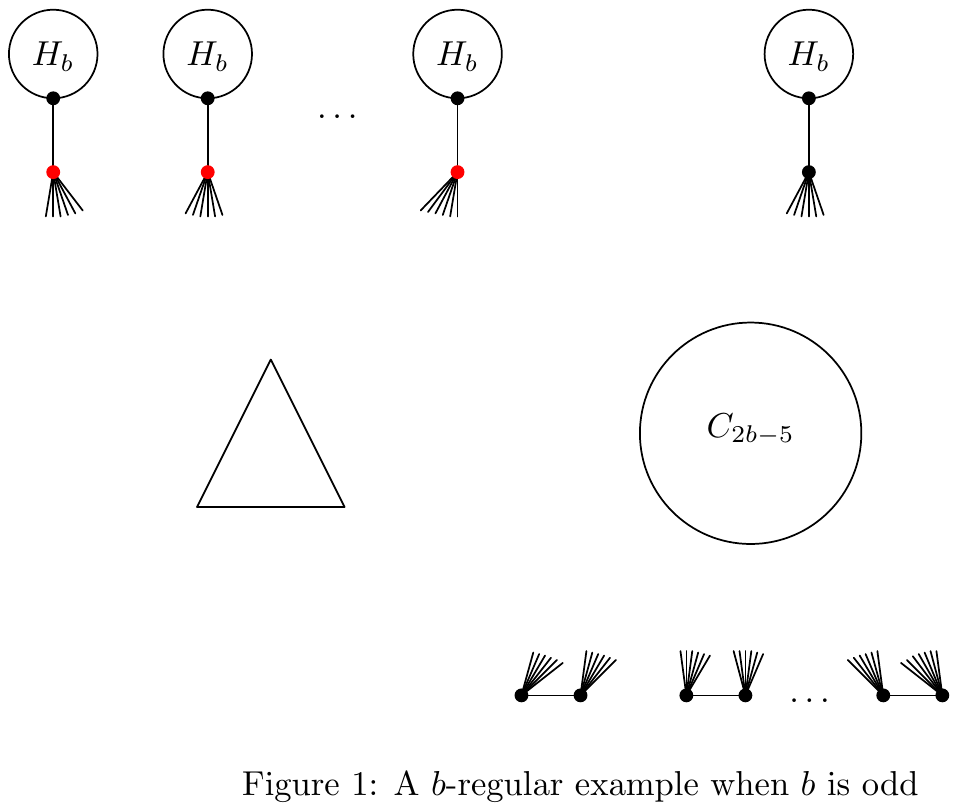}
\vspace{-360pt}
\end{figure}
The obtained graph is $b$-regular and connected, and $X$ is a Godsil-McKay switching set.
Furthermore, by deleting the $b-2$ vertices of $W$, the corresponding $b-2$ copies of $H_b$, 
the triangle and the remainder become $b$ components, each having an odd number of vertices.
Consequently %by Tutte's theorem
this graph does not have a perfect matching.
However, after performing a Godsil-McKay switch, one easily finds (many) perfect matchings.
This concludes the proof for the odd case.

When $b$ is even, we make a few small changes.
First, the graphs $\widetilde{H_b}$ should be replaced with a graph obtained as follows.
Delete one edge from the complete graph $K_{b+1}$, add an additional vertex $v$ and connect it to the two vertices
of degree $b-1$.
Then, since $b+1$ is odd, any perfect matching must contain one of these two new edges,
which precludes any additional edge incident with $v$ from being in a matching, just as before.
The switching set $X$ gets $2b-4$ vertices and induces $K_3+C_{2b-7}$.

Lastly, we make a small alteration to the final step, where we increased the degrees of the vertices in
the larger cycle to $b$.
Because $b$ is even, we can do this without adding an additional $\widetilde{H_b}$.
Instead, we add a cycle $C_{b-2}$ to $Y$ and join each vertex in the added cycle to $b-2$ vertices of the
larger cycle in $X$, such that all degrees become $b$.
With the mentioned modifications we complete the proof of the theorem by imitating the above steps for the even case.
\hfill$\Box$
\\[5pt]
We remark that in some cases we can do better by taking a regular graph of larger degree on the switching set $X$.
However, we do not know how to modify the construction to make it work for $b\leq 4$.
In fact, if $b\leq 2$ there exist no nonisomorphic cospectral graphs (see~\cite{vDH}),
and if $b=3$ it can be seen that there cannot exist a Godsil-McKay switch between a graph with a
perfect matching and one with none.
Moreover, Stephen Hartke checked by computer all $3$-regular pairs of cospectral graphs
on at most 20 vertices and found no example.
So it is not unlikely that for degree $3$ there doesn't exist such a pair of graphs.
\\[5pt]
{\bf \large Acknowledgement.}
The result of this note was found at the sixth Eml\'ekt\'abla workshop that took place in July 2014 in Hungary.
The authors are grateful to the organizers of this meeting for creating an excellent atmosphere for problem solving.

\end{document}